\documentclass[a4paper,12pt]{article}
\usepackage{comment}
\usepackage{cite}
\usepackage{amsmath}
\usepackage{amssymb}
\usepackage{amsfonts}
\usepackage[T1]{fontenc}
\usepackage[utf8]{inputenc}
\usepackage{graphicx}
\usepackage{fancyhdr}
\usepackage{float}
\usepackage{xcolor}
\usepackage{authblk}
\usepackage{mathrsfs}
\usepackage{empheq}
\usepackage[hyphens]{url}
\usepackage{hyperref} 
\usepackage[]{breakurl}
\usepackage{amsmath}
\usepackage{amssymb}

\usepackage{geometry}
\begin{document}

\title{On an upper bound for central binomial coefficients and Catalan numbers}

\author[$\dagger$]{Jean-Christophe {\sc Pain}$^{1,2,}$\footnote{jean-christophe.pain@cea.fr}\\
\small
$^1$CEA, DAM, DIF, F-91297 Arpajon, France\\
$^2$Universit\'e Paris-Saclay, CEA, Laboratoire Mati\`ere en Conditions Extr\^emes,\\ 
F-91680 Bruy\`eres-le-Ch\^atel, France
}

\date{}

\maketitle

\begin{abstract}
Recently, Agievich proposed an interesting upper bound on binomial coefficients in the de Moivre-Laplace form. In this article, we show that the latter bound, in the specific case of a central binomial coefficient, is larger than the one proposed by Sasvari and obtained using the Binet formula for the Gamma function. In addition, we provide the expression of the next-order bound and apply it to Catalan numbers $C_n$. The bounds are very close to the exact value, the difference decreasing with $n$ and with the order of the upper bound. 
\end{abstract}

\section{Introduction}

Agievich recently published the upper bound \cite{Agievich2022}:
\begin{equation}\label{agie}
    \binom{n}{k}\leq\frac{2^n}{\sqrt{\pi n/2}}\,\exp\left[-\frac{2}{n}\left(k-\frac{n}{2}\right)^2+\frac{23}{18n}\right],
\end{equation}
where $(n,k)\in\mathbb{N}$ and $0\leq k\leq n$. One has, as a special case of Eq. (\ref{agie}), the following bound:
\begin{equation}
    \binom{2n}{n+k}\leq\frac{2^{2n}}{\sqrt{\pi n}}\,\exp\left(-\frac{k^2}{n}+\frac{23}{36n}\right),
\end{equation}
yielding, for $k=0$: 
\begin{equation}
    \binom{2n}{n}\leq\frac{2^{2n}}{\sqrt{\pi n}}\,\exp\left(\frac{23}{36n}\right),
\end{equation}
and thus, for Catalan numbers \cite{Stanley2015}:
\begin{equation}
C_n\leq\frac{2^{2n}}{(n+1)\sqrt{\pi n}}\,\exp\left(\frac{23}{36n}\right).
\end{equation}
In 1999, Sasvari obtained \cite{Sasvari1999}, using the Binet integral representation of the Gamma function \cite{Sasvari1999b}:
\begin{equation}
    \frac{4^n}{\sqrt{\pi n}}\,\exp\left(-\frac{1}{8n}\right)\leq\binom{2n}{n}\leq\frac{4^n}{\sqrt{\pi n}}\,\exp\left(-\frac{1}{8n}+\frac{1}{192n^3}\right).
\end{equation}
It appears that Sasvari's upper bound 
\begin{equation}
    U_S(n)=\frac{4^n}{\sqrt{\pi n}}\,\exp\left(-\frac{1}{8n}+\frac{1}{192n^3}\right)
\end{equation}
is tighter than Agievich's one 
\begin{equation}
    U_A(n)=\frac{2^{2n}}{\sqrt{\pi n}}\,\exp\left[\frac{23}{36n}\right].
\end{equation}
This can be checked easily by considering (for $x$ real strictly positive) the function
\begin{equation}
    f(x)=\frac{23}{36x}+\frac{1}{8x}-\frac{1}{192x^3}=\frac{55}{72x}-\frac{1}{192x^3}.
\end{equation}
We find that $f$ is decreasing for
\begin{equation}
    x>\left(\frac{9}{440}\right)^{1/2}=\frac{3}{2\,\sqrt{110}}<1
\end{equation}
and since $\lim_{x\rightarrow\infty} f(x)=0$, we have necessarily $f(x)>0$ and therefore $U_S(n)<U_A(n)$. The values of $U_A(n)$ and $U_S(n)$ for the ten first values of $n$ are displayed in Table \ref{tab0}.

\begin{table}
\begin{center}
\begin{tabular}{cccc}\hline\hline
$n$ & $\binom{2n}{n}$ & $U_A(n)$ & $U_S(n)$ \\\hline\hline
1 & 2 & 4.275146228 & 2.001982123 \\
2 & 6 & 8.785429701 & 6.000250574 \\
3 & 20 & 25.79485257 & 20.00011914 \\
4 & 70 & 84.72303658 & 70.00010214 \\
5 & 252 & 93.5845534 & 252.0001224 \\
6 & 924 & 1049.430558 & 924.0001819 \\
7 & 3432 & 3827.665444 & 3432.000314 \\
8 & 12870 & 14159.34751 & 12870.00061 \\
9 & 48620 & 52926.51245 & 48620.00127 \\
10 & 184756 & 199421.3118 & 184756.0029 \\\hline\hline
\end{tabular}
\caption{Values of Agierich's bound $U_A(n)$ (column 2) and Sasvari's bound $U_S(n)$ (column 3) with 10 significant digits, for the first 10 values of $n$. Comparison with the exact values of $\binom{2n}{n}$ (first column).}\label{tab0}
\end{center}
\end{table}

\section{Next-order bound}

One has the general upper bound (note that here is a typographical error in Ref. \cite{Sasvari1999} since $d_r^{s}$ should be $d_r^{2s}$):
\begin{equation}
\binom{rs}{s}<c_r\,d_r^{2s}\frac{1}{\sqrt{s}}\,\exp\left[D_{2N}(s,r)\right],
\end{equation}
with
\begin{equation}
D_N(s,r)=\sum_{j=1}^N\frac{B_{2j}}{2j(2j-1)}\left[\frac{1}{(rs)^{2j-1}}-\frac{1}{(s)^{2j-1}}-\frac{1}{[(r-1)s]^{2j-1}}\right]
\end{equation}
where $B_k$ are the Bernoulli numbers defined by
\begin{equation}
\frac{t}{e^t-1}=1-\frac{t}{2}+\sum_{j=1}^{\infty}\frac{B_{2j}}{(2j)!}t^{2j}
\end{equation}
and
\begin{equation}
c_r=\frac{1}{\sqrt{2\pi\left(1-\frac{1}{r}\right)}}
\end{equation}
as well as
\begin{equation}
d_r=\frac{r-1}{\left(1-\frac{1}{r}\right)}.
\end{equation}

Setting $s=n$ and $r=2$, one gets
\begin{equation}
D_N(n,2)=\sum_{j=1}^N\frac{B_{2j}}{2j(2j-1)}\left[\frac{1}{(2n)^{2j-1}}-\frac{2}{n^{2j-1}}\right],
\end{equation}
where
$c_2=1/\sqrt{\pi}$ and $d_2=2$. One has
\begin{equation}
D_{2N}(n,2)=\sum_{j=1}^{2N}\frac{B_{2j}}{j(2j-1)}\left[\frac{1}{2^{2j}-1}\right]\frac{1}{n^{2j-1}}
\end{equation}
and thus
\begin{equation}
\binom{2n}{n}<\frac{2^n}{\sqrt{\pi n}}\,\exp\left\{\sum_{j=1}^N\frac{B_{2j}}{j(2j-1)}\left[\frac{1}{2^{2j}-1}\right]\frac{1}{n^{2j-1}}\right\}.
\end{equation}
One has in particular
\begin{equation}
D_2(n,2)=-\frac{3B_2}{4n}-\frac{5B_4}{32n^3}
\end{equation}
and for the next order
\begin{equation}
D_4(n,2)=-\frac{3B_2}{4n}-\frac{5B_4}{32n^3}-\frac{21B_6}{320n^5}-\frac{255B_8}{7168n^7},
\end{equation}
where $B_2=1/6$, $B_4=-1/30$, $B_6=-1/42$ and $B_8=-1/30$. This yields
\begin{equation}
D_2(n,2)=-\frac{1}{8n}+\frac{1}{192n^3}
\end{equation}
and for the next order
\begin{equation}
D_4(n,2)=-\frac{1}{8n}+\frac{1}{192n^3}-\frac{1}{640n^5}+\frac{17}{14336n^7},
\end{equation}

Values of $\binom{2n}{n}\,\sqrt{\pi n}/2^{2n}$ and upper bounds $\exp[D_2(n,2)]$ and $\exp[D_4(n,2)]$ are displayed in Table \ref{tab1}. Similarly, one has, for Catalan numbers, the following upper bound 
\begin{equation}\label{bou1}
    C_n<\frac{2^{2n}}{(n+1)\,\sqrt{\pi n}}\,\exp\left[-\frac{1}{8n}+\frac{1}{192n^3}\right]
\end{equation}
and for the next order
\begin{equation}\label{bou2}
    C_n<\frac{2^{2n}}{(n+1)\,\sqrt{\pi n}}\,\exp\left[-\frac{1}{8n}+\frac{1}{192n^3}-\frac{1}{640n^5}+\frac{17}{14336n^7}\right];
\end{equation}
both are given in Table \ref{tab2} (third and fourth columns respectively).

\begin{table}
\begin{center}
\begin{tabular}{cccc}\hline\hline
$n$ & Exact & Bound $N=1$ & Bound $N=2$ \\\hline\hline
1 & 0.88622692545276 & 0.88710523105688 & 0.88677114441088 \\
2 & 0.93998560298663 & 0.94002485899037 & 0.93998766871395 \\
3 & 0.95936878869983 & 0.95937450378689 & 0.95936885517397 \\
4 & 0.96931069971395 & 0.96931211408847 & 0.96931070519255 \\
5 & 0.97535007714523 & 0.97535055078797 & 0.97535007791724 \\
6 & 0.97940560431422 & 0.97940579711995 & 0.97940560446817 \\
7 & 0.98231617716265 & 0.98231626711057 & 0.98231617720181 \\
8 & 0.98450640547183 & 0.98450645187199 & 0.98450640548375 \\
9 & 0.98621413686019 & 0.98621416271614 & 0.98621413686436 \\
10 & 0.98758292882616 & 0.98758294414165 & 0.98758292882778 \\\hline\hline
\end{tabular}
\caption{Values of $\binom{2n}{n}\,\sqrt{\pi n}/2^{2n}$ and upper bounds $\exp\left[D_2(n,2)\right]$ and $\exp\left[D_4(n,2)\right]$ with 14 significant digits.}\label{tab1}
\end{center}
\end{table}

\begin{table}
\begin{center}
\begin{tabular}{cccc}\hline\hline
$n$ & Exact & Bound $N=1$ & Bound $N=2$ \\\hline\hline
$C_1$ & 1 & 1.0009910617460 & 1.0006140853347 \\
$C_2$ & 2 & 2.0000835246915 & 2.0000043952318 \\
$C_3$ & 5 & 5.0000297856629 & 5.0000003464472 \\
$C_4$ & 14 & 14.000020428169 & 14.000000079129 \\
$C_5$ & 42 & 42.000020395749 & 42.000000033244 \\
$C_6$ & 132 & 132.00002598551 & 132.00000002075 \\
$C_7$ & 429 & 429.00003928232 & 429.00000001710 \\
$C_8$ & 1430 & 1430.0000673964 & 1430.0000000173 \\
$C_9$ & 4862 & 4862.0001274689 & 4862.0000000205 \\
$C_{10}$ & 16796 & 16796.000260473 & 16796.000000028 \\\hline\hline
\end{tabular}
\caption{The ten first Catalan numbers (second column) and upper bounds given by (\ref{bou1}) and (\ref{bou2}) (columns 3 and 4 respectively) with 14 significant digits.}\label{tab2}
\end{center}
\end{table}

It is worth mentioning that nice bounds were also recently published by St$\mathrm{\breve{a}}$nic$\mathrm{\breve{a}}$ for specific binomial coefficients \cite{Stanica2001}.

\section{Conclusion}

Agievich proposed a clever derivation of an upper bound on binomial coefficients in the de Moivre-Laplace form. It appears that, in the particular case of a central binomial coefficient, the corresponding upper bound is larger than the one published by Sasvari in 1999. We gave the expression of the next-order bound in Sasvari's approach and applied it to Catalan numbers. The bounds are very close to the exact value, especially for high integer orders of the Catalan numbers. The next term in the expansion of the upper bound reduces drastically the difference with the exact value.

\end{document}